\documentclass[a4paper,12pt]{amsart}

\makeatletter

\@setfontsize\normalsize\@xiipt{15.5}%

\renewcommand{\theenumi}{\alph{enumi}}

\renewcommand{\theenumii}{\roman{enumii}}

\renewcommand{\p@enumii}{\theenumi.}
\renewcommand{\theenumiii}{\Alph{enumiii}}

\renewcommand{\p@enumiii}{\theenumi.\theenumii.}

\renewcommand{\p@enumiv}{\p@enumiii\theenumiii.}

\@addtoreset{equation}{section}

\RequirePackage{amsthm}

\renewenvironment{proof}[1][\proofname]{\par
  \normalfont
  \topsep6\p@\@plus6\p@ \trivlist
  \item[\hskip\labelsep\slshape
    #1\@addpunct{.}]\ignorespaces
}{%
  \qed\endtrivlist
}

\theoremstyle{plain}

\newtheorem{theorem}{Theorem}[section]
\newtheorem{proposition}[theorem]{Proposition}
\newtheorem{lemma}[theorem]{Lemma}
\newtheorem{corollary}[theorem]{Corollary}

\theoremstyle{definition}

\newtheorem{definition}[theorem]{Definition}

\newtheorem{remark}[theorem]{Remark}

\RequirePackage{amsmath}
\RequirePackage{amssymb}

\RequirePackage{xypic}

\renewcommand{\emptyset}{\varnothing}

\newcommand{\Union}{\bigcup\limits}

\newcommand{\Inter}{\bigcap\limits}

\newcommand{\C}{\mathbb{C}}

\newcommand{\R}{\mathbb{R}}

\newcommand{\Z}{\mathbb{Z}}

\newcommand{\BDC}{\mathbf{D}^{\mathrm{b}}}

\newcommand{\DSum}{\bigoplus}

\newcommand{\ilim}[1][]{\mathop{\varinjlim}\limits_{#1}}

\renewcommand{\to}[1][]{\xrightarrow[#1]{}}
\newcommand{\from}[1][]{\xleftarrow[#1]{}}
\newcommand{\isoto}[1][]{\xrightarrow[#1]{\sim}}

\newcommand{\Endo}[1][]{\mathrm{End}_{\raise1.5ex\hbox to.1em{}#1}}

\newcommand{\Hom}[1][]{\mathrm{Hom}_{\raise1.5ex\hbox to.1em{}#1}}

\newcommand{\RHom}[1][]{\mathrm{RHom}_{\raise1.5ex\hbox to.1em{}#1}}

\newcommand{\Ext}[2][]{\mathrm{Ext}_{\raise1.5ex\hbox to.1em{}#1}^{#2}}

\newcommand{\THom}[1][]{\mathrm{THom}_{\raise1.5ex\hbox to.1em{}#1}}

\newcommand{\Tens}[1][]{\mathbin{\otimes_{\raise1.5ex\hbox to-.1em{}#1}}}

\newcommand{\LTens}[1][]{\mathbin{\otimes_{\raise1.5ex\hbox to-.1em{}#1}^{L}}}

\newcommand{\Tor}[2][]{\mathrm{Tor}^{\raise1.5ex\hbox to.1em{}#1}_{#2}}

\def\shb{\mathcal{B}}

\def\shm{\mathcal{M}}

\newcommand{\sect}{\Gamma}
\newcommand{\rsect}{\mathrm{R}\Gamma}

\renewcommand{\hom}[1][]{{\mathcal{H}om}_{\raise1.5ex\hbox to.1em{}#1}}

\newcommand{\rhom}[1][]{{R\mathcal{H}om}_{\raise1.5ex\hbox to.1em{}#1}}

\newcommand{\ext}[2][]{{\mathcal{E}xt}_{\raise1.5ex\hbox to.1em{}#1}^{#2}}

\newcommand{\thom}[1][]{{T\mathcal{H}om}_{\raise1.5ex\hbox to.1em{}#1}}

\newcommand{\tens}[1][]{\mathbin{\otimes_{\raise1.5ex\hbox to-.1em{}#1}}}

\newcommand{\ltens}[1][]{\mathbin{\otimes_{\raise1.5ex\hbox to-.1em{}#1}^{L}}}

\newcommand{\tor}[2][]{{\mathcal{T}or}^{\raise1.5ex\hbox to.1em{}#1}_{#2}}

\DeclareMathOperator{\supp}{supp}

\newcommand{\roim}[1]{{R#1}_*}
\newcommand{\reim}[1]{{R#1}_!}

\newcommand{\opb}[1]{#1^{-1}}

\newcommand{\epb}[1]{#1^{!}}

\DeclareMathOperator{\ori}{or}

\newcommand{\GHom}[1][]{\mathrm{GHom}_{\raise1.5ex\hbox to.1em{}#1}}

\newcommand{\GExt}[2][]{\mathrm{GExt}_{\raise1.5ex\hbox to.1em{}#1}^{#2}}

\newcommand{\FHom}[1][]{\mathrm{FHom}_{\raise1.5ex\hbox to.1em{}#1}}

\newcommand{\ghom}[1][]{{\mathcal{GH}om}_{\raise1.5ex\hbox to.1em{}#1}}

\newcommand{\gext}[2][]{{\mathcal{GE}xt}_{\raise1.5ex\hbox to.1em{}#1}^{#2}}

\newcommand{\fhom}[1][]{{\mathcal{FH}om}_{\raise1.5ex\hbox to.1em{}#1}}

\newcommand{\tenstop}[1][]{\mathbin{\hat{\otimes}_{\raise1.5ex\hbox to-.1em{}#1}}}

\newcommand{\homtop}[1][]{\mathcal{L}_{\raise1.5ex\hbox to.1em{}#1}}

\newcommand{\Homtop}[1][]{\mathrm{L}_{\raise1.5ex\hbox to.1em{}#1}}

\newcommand{\D}{\mathcal{D}}

\renewcommand{\O}{\mathcal{O}}

\newcommand{\Db}{\mathcal{D}b}

\DeclareMathOperator{\chv}{char}

\def\absdoim#1{\underline{#1}_*}
\def\reldoim[#1]#2{\underline{#2}_{|{#1}*}}
\def\doim{\@ifnextchar [{\reldoim}{\absdoim}}

\def\absdeim#1{\underline{#1}_*}
\def\reldeim[#1]#2{\underline{#2}_{|{#1}*}}
\def\deim{\@ifnextchar [{\reldeim}{\absdeim}}

\def\absdopb#1{\underline{#1}^{-1}}
\def\reldopb[#1]#2{\underline{#2}_{|{#1}}^{-1}}
\def\dopb{\@ifnextchar [{\reldopb}{\absdopb}}

\def\absboim#1{\underline{\underline{#1}}_*}
\def\relboim[#1]#2{\underline{\underline{#2}}_{|{#1}*}}
\def\boim{\@ifnextchar [{\relboim}{\absboim}}

\def\absbeim#1{\underline{\underline{#1}}_*}
\def\relbeim[#1]#2{\underline{\underline{#2}}_{|{#1}*}}
\def\beim{\@ifnextchar [{\relbeim}{\absbeim}}

\def\absbopb#1{\underline{\underline{#1}}^*}
\def\relbopb[#1]#2{\underline{\underline{#2}}_{|{#1}}^*}
\def\bopb{\@ifnextchar [{\relbopb}{\absbopb}}

\makeatother

\newcommand{\ft}{{{}^t\!f'}}

\newcommand{\fp}{{f_\pi}}

\newcommand{\inv}{{-1}}

\newcommand{\fut}[1]{{#1}^{\uparrow}}

\newcommand{\pas}[1]{{#1}^{\downarrow}}

\newcommand{\cone}{\gamma}

\newcommand{\dcone}{\lambda}

\newcommand{\Vcone}{C}

\newcommand{\Int}{\operatorname{Int}}

\date{}

\address{Institut de Math\'ematiques;
Analyse Alg\'ebrique (Case 82);
Universit\'e Pierre et Marie Curie;
4, place Jussieu;
F-75252 Paris Cedex 05}

\email{dagnolo@math.jussieu.fr, schapira@math.jussieu.fr}

\title{Global propagation on causal manifolds}

\author{Andrea D'Agnolo \and Pierre Schapira}

\dedicatory{Dedicated to Professor Mikio Sato on the occasion of his 
70th birthday}

\begin{document}

\maketitle

\section*{Introduction}

The micro-support of sheaves (see~\cite{Kashiwara-Schapira90}) is a
tool to describe local propagation results.  A natural problem is then
to give sufficient conditions to get global propagation results from the
knowledge of the micro-support.  This is the aim of this paper.

A propagator on a real manifold $M$ is the data of a pair
$(Z,\dcone)$, where $Z\subset M\times M$ is a closed subset containing
the diagonal, $\dcone$ is a closed cone of the cotangent bundle to
$M$, and some relation holds between $\dcone$ and the micro-support of
the constant sheaf along $Z$.  In this framework, we prove that if $F$
is a sheaf on $M$ whose micro-support does not intersect $-\dcone$
outside of the zero-section, then the restriction morphism
$\rsect(M;F) \to \rsect(U;F)$ is an isomorphism, as soon as $M
\setminus U$ is $Z$-proper.  This last condition means that the
forward set $\fut D = \{y\in M\colon (x, y) \in Z\text{ for some }x\in
D\}$ of any compact set $D\subset M$ should intersect $M \setminus U$
in a compact set, and the backward set $\pas{(M \setminus U)} = \{x\in
M\colon (x, y) \in Z\text{ for some }y\notin U\}$ should not contain
any connected component of $M$.

As an application, we consider the problem of global existence for
solutions to hyperbolic systems (in the hyperfunction and distribution
case), along the lines of Leray~\cite{Leray55}.  Causal manifolds, and
in particular homogeneous causal manifolds as considered by Faraut et
al., give examples of manifolds to which our results apply.

\section{Statement of the results}

\subsection{Normal cones}

A subset $\Vcone$ of a finite dimensional real vector space $V$ is 
called a cone (or a conic subset), if $\R^+ \cdot \Vcone\subset 
\Vcone$.  A cone $\Vcone\subset V$ is called convex if 
$\Vcone+\Vcone\subset \Vcone$, and proper if $\overline\Vcone\cap 
-\overline\Vcone \subset \{0\}$.  We also use the notation 
$\Vcone^a=-\Vcone$.  Denoting by $V^*$ the dual of $V$, the polar to a 
cone $\Vcone\subset V$ is the conic subset of $V^*$ defined by 
$\Vcone^\circ=\{\xi\colon \langle\xi,v\rangle\geq 0\text{ for every 
}v\in \Vcone\}$.  One checks that $(C^\circ)^\circ$ is the closure of 
the convex envelop to $C$, and that the polar to a proper convex cone 
is a closed proper convex cone.

Let $M$ be a $C^\infty$-manifold.  If $q\colon E\to M$ is a vector
bundle, one naturally extends the above notions to subsets of $E$.
For example, $\cone\subset E$ is a cone if $\cone_{x}:=\cone\cap
q^\inv(x)$ is a cone in $E_{x}$ for any $x\in X$.  We identify $M$ to
the zero-section of $q$, and for $\cone \subset E$ we set $\dot\cone =
\cone \setminus M$.

Denote by $\tau\colon TM\to M$ and $\pi \colon T^*M \to M$ the tangent
and cotangent bundle to $M$, respectively.
Following~\cite[Definition~4.1.1]{Kashiwara-Schapira90}, $C(A,B)$
denotes the Whitney normal cone of $A,B\subset M$, a closed cone of
$TM$.  Recall that if $(x)$ is a local coordinate system in $M$, then
$(x_{\circ};v_{\circ})\in C(A,B)$ if and only if there exists a
sequence $(a_{n},b_{n},c_{n})$ in $A\times B\times \R^+$ such that
$$
a_{n} \to x_{\circ}, \quad b_{n} \to x_{\circ}, \quad c_{n}(a_{n} -
b_{n}) \to v_{\circ}.
$$
If $N\subset M$ is a smooth submanifold, $C_{N}(A)$ is the projection
of $C(N,A)$ in $T_{N}M$, the normal bundle to $N$ in $M$.

The {\em strict normal cone} to $A\subset M$ is defined
in~\cite[Definition~5.3.6]{Kashiwara-Schapira90} by $N(A)=TM\setminus
C(M\setminus A,A)$.  Recall that if $(x)$ is a local coordinate system
in $M$, then $(x_{\circ};v_{\circ})\in N(A)$ if and only if there
exists an open cone $\Vcone$ containing $v_{\circ}$, and a
neighborhood $U$ of $x_{\circ}$, such that
\begin{equation}
    U \cap \bigl( (A\cap U)+\Vcone \bigr) \subset A.
    \label{eq:strict}
\end{equation}
Note that $N(A)$ is an open convex cone of $TM$, $N(M\setminus A) =
N(A)^{a}$, and $N_{x}(A) \neq T_{x}M$ if and only if $x$ is in the
topological boundary of $A$.

\subsection{Micro-support}

Let $M$ be a $C^\infty$-manifold.  Let $k$ be a field, and denote by
$\BDC(k_{M})$ the bounded derived category of sheaves of $k$-vector
spaces on $M$.  Following~\cite[Chapter~5]{Kashiwara-Schapira90}, to
$F\in\BDC(k_{M})$ one associates its {\em micro-support} $SS(F)$, a
closed conic involutive subset of $T^*M$.  Recall that $T^*M\setminus
SS(F)$ describes the (co)directions of propagation for the cohomology
of $F$, stable by small perturbations.  More precisely, $p\notin
SS(F)$ if and only if there exists an open neighborhood $\Omega$ of
$p$ such that for any $x\in\pi(\Omega)$ and any $C^\infty$-function
$\varphi$ on $M$ with $\varphi(x)=0$, $d\varphi(x)\in\Omega$, one has
\begin{equation}
    \left( \rsect_{\{\varphi\geq 0\}}F \right)_{x}=0,
    \label{eq:SS}
\end{equation}
where $\rsect_W$ denotes the derived functor of sections with support
on a closed subset $W\subset M$, and we write for short $\{\varphi\geq
0\} = \{y\in M\colon \varphi(y)\geq 0\}$.  This is indeed a
propagation requirement, since the above vanishing can be restated by
asking that the natural restriction morphism
$$
\ilim[U\owns x] H^{j}(U;F) \to \ilim[U\owns x] H^{j}(U\cap\{\varphi <
0\};F)
$$
is an isomorphism for any $j\in\Z$.  This implies that ``sections'' of
$F$ on $U\cap\{\varphi < 0\}$ extend to a neighborhood of $x$.

If $A\subset M$ is a locally closed subset, denote by $k_{A}$ the
sheaf on $M$ which is zero on $M\setminus A$, and constant with fiber
$k$ on $A$.  Recall that if $U\subset M$ is an open subset, and
$W\subset M$ is a closed subset, one has the estimates:
\begin{equation}
    SS(k_{U})\subset N(U)^{\circ a}, \qquad SS(k_{W})\subset
    N(W)^{\circ}.
    \label{eq:SSinN}
\end{equation}

\subsection{Propagators}

Let $M$ be a $C^\infty$-manifold.  Denote by $\Delta\subset M\times M$
the diagonal, and by $q_{1}$ and $q_{2}$ the first and second
projection from $M\times M$ to $M$.

\begin{definition}
    Let $Z\subset M\times M$ be a closed subset.  We say that a
    locally closed subset $A\subset M$ is {\em $Z$-proper} if
    \begin{itemize}
	\item [(i)] $q_{1}$ is proper on $Z\cap q_{2}^\inv (\overline
	A)$,

	\item [(ii)] $q_{1} \bigl( Z\cap q_{2}^\inv(\overline A)
	\bigr)$ does not contain any connected component of $M$.
    \end{itemize}
\end{definition}

Given $Z\subset M\times M$ as above, to a subset $A\subset M$, we
associate
\begin{eqnarray*}
    \pas A &=& q_{1}(Z\cap q_{2}^\inv A), \\
    \fut A &=& q_{2}(Z\cap q_{1}^\inv A),
\end{eqnarray*}
and we set $\pas x = \pas{\{x\}}$, $\fut x = \fut {\{x\}}$.  With
these notations, a subset $A\subset M$ is $Z$-proper if and only if:
(i) $\fut D \cap \overline A$ is compact for any compact subset $D$ of
$M$, (ii) $\pas{\overline A}$ does not contain any connected component
of $M$.

\begin{definition}
    Let $Z$ be a closed subset of $M\times M$, and $\dcone$ a closed
    cone of $T^*M$.  We say that the pair $(Z,\dcone)$ is a {\em
    propagator} on $M$ if
    \begin{itemize}
	\refstepcounter{equation}\label{hy:Zdelta} \item
	[\eqref{hy:Zdelta}] $\Delta\subset Z$,

	\refstepcounter{equation}\label{hy:ZsubGG} \item
	[\eqref{hy:ZsubGG}] $SS(k_{Z})\subset T^*M\times\dcone$,

	\refstepcounter{equation}\label{hy:ZtransM} \item
	[\eqref{hy:ZtransM}] $SS(k_{Z}) \cap (T^{*}M\times M) \subset
	M\times M$,

	\refstepcounter{equation}\label{hy:ZtransN} \item
	[\eqref{hy:ZtransN}] $SS(k_{Z}) \cap (M\times T^{*}M) \subset
	M\times M$.
    \end{itemize}
    (As for \eqref{hy:ZtransM} and \eqref{hy:ZtransN}, recall that we
    identify the zero-section of $T^{*}M$ to $M$.)  We say that
    $(Z,\dcone)$ is a {\em convex propagator} on $M$ if it is a
    propagator and moreover
    \begin{itemize}
	\refstepcounter{equation}\label{hy:Gconvex} \item
	[\eqref{hy:Gconvex}] $\dcone$ is a proper convex cone.
    \end{itemize}
\end{definition}

\subsection{Propagation theorems}

\begin{theorem}
    \label{th:propag0}
    Let $(Z,\dcone)$ be a propagator on $M$.  Let $F\in\BDC(k_{M})$,
    and assume that $\supp(F)$ is $Z$-proper and
    $SS(F)\cap\dcone^{a}\subset M$.  Then $$\rsect(M;F)=0.$$
\end{theorem}

Part~(i) of the following corollary partially extends to manifolds
Proposition~5.2.1 of~\cite{Kashiwara-Schapira90} which only considered
an affine situation, with $\dcone$ constant along the fibers.  (See
Remark~\ref{re:KS} for further comments.)

\begin{corollary}
    \label{co:propag}
    Let $(Z,\dcone)$ be a convex propagator on $M$.  Let
    $F\in\BDC(k_{M})$, and assume that $SS(F)\cap\dcone^{a}\subset M$.
    \begin{itemize}
	\item [(i)] Let $W$ be a closed subset of $M$ which is
	$Z$-proper and satisfies $SS(k_{W})\subset\dcone^{a}$.  Then
	$$\rsect_{W}(M;F)=0.$$

	\item [(ii)] Let $U$ be an open subset of $M$ which is
	$Z$-proper and satisfies $SS(k_{U})\subset\dcone$.  Then
	$$\rsect(M;F_{U})=0.$$
    \end{itemize}
\end{corollary}

Note that (i) and (ii) are equivalent to
\begin{eqnarray*}
    \rsect(M;F) &\isoto& \rsect(M\setminus W;F), \\
    \rsect(M;F) &\isoto& \rsect(M\setminus U;F),
\end{eqnarray*}
respectively.  In other words, ``sections'' of $F$ on $M\setminus W$
(or on a neighborhood of $M\setminus U$) extend uniquely to $M$.

The following result deals with the case where $\dcone$ is not convex,
but is covered by a finite union of convex cones.  A situation which
appears for example in dealing with the Cauchy problem, real or
complex.

\begin{corollary}
    \label{co:PropagFinite}
    Let $I$ be a finite set.  For $j\in I$, let $U_{j}$ be an open
    subset of $M$, and set $N=M\setminus\Union_{j\in I}U_{j}$.  For
    any $J \subset I$, $J \neq \emptyset$, let $(Z_{J},\dcone_{J})$ be
    a convex propagator, and set $U_{J}=\Inter_{j\in J}U_{j}$.  Let
    $F\in\BDC(k_{M})$.  Assume
    \begin{itemize}
	\item [(i)] $U_{J}$ is $Z_{J}$-proper, and $SS(k_{U_{J}})
	\subset \dcone_{J}$,

	\item [(ii)] $SS(F)\cap\dcone_{J}^{a} \subset M$.
    \end{itemize}
    Then, one has the isomorphism $$ \rsect(M;F) \isoto \rsect(N;F).
    $$
\end{corollary}

\begin{remark}
    \label{re:KS}
    Theorem~\ref{th:propag0} does not allow one to recover
    Proposition~5.2.1 of~\cite{Kashiwara-Schapira90}, since our
    hypotheses are stronger.  More precisely, we require $\dcone$
    closed proper convex and $SS(F)\cap\dcone^{a}\subset M$, while in
    loc.~cit.\ one only assumes $\dcone=\cone^\circ$, for $\cone$
    closed proper convex in $TM$, and $SS(F)\cap \Int(\dcone^{a}) =
    \emptyset$.  Let us give an example which shows that, in general,
    it is not possible to replace the hypothesis
    $SS(F)\cap\dcone^{a}\subset M$ by the hypothesis $SS(F)\cap
    \Int(\dcone^{a}) = \emptyset$.

    Let $M=\R\times S^1$ be an infinite cylinder.  Using the
    identification $T^*M = M\times (\R\times\R)$, set $Z =
    \{(x_{1},\theta_{1},x_{2},\theta_{2})\in M\times M\colon x_{1} =
    x_{2}\}$, $\dcone = \{(x,\theta; \xi, \tau)\in T^*M \colon
    \tau\geq 0\}$, $W=\{0\}\times S^1\subset M$.  Clearly,
    $(Z,\dcone)$ is a propagator on $M$, and $W$ is $Z$-proper.
    Setting $F=k_{W}$, one has $SS(F) = \{(x,\theta; \xi, \tau)\in
    T^*M \colon x = \tau = 0\}$, and hence
    \begin{align*}
	& SS(F) \cap \Int(\dcone^{a}) = \emptyset, \\
	& SS(F) \cap \dcone^a \not\subset M, \\
	& \sect(M;F) \neq 0.
    \end{align*}
\end{remark}

\section{Proof of the results}

\subsection{Review on sheaves}

Let $f:N\to M$ be a morphism of $C^\infty$ manifolds.  We will
consider the usual operations $\roim f$, $\reim f$, $\opb f$, $\epb
f$, $\tens$, $\rhom$ of sheaf theory.  If $F\in\BDC(k_{M})$, we set
$D'F=\rhom(F,k_{M})$.  We also make use of the absolute and relative
dualizing complexes denoted $\omega_{M}$ and $\omega_{N/M}$,
respectively.  Recall that if $f$ is smooth, then
$\omega_{N/M}=\ori_{N/M}[\dim N - \dim M]$, where $\ori_{N/M}$ denotes
the relative orientation sheaf.

We will need the following lemma.

\begin{lemma}
    \label{le:summand}
    Let $Z$ be a closed subset of $M\times M$ containing the diagonal
    $\Delta$.  Then
    \begin{itemize}
	\item [(i)] $k_{M}$ is a direct summand of
	$\roim{q_{2}}k_{Z}$,

	\item [(ii)] $\omega_{M}$ is a direct summand of
	$\reim{q_{2}}\omega_{Z}$.
    \end{itemize}
\end{lemma}

\begin{proof}
    Since the arguments are similar, we will prove only (ii).  Set
    $\widetilde q_{2}=q_{2}\vert_{Z}$, and denote by $i\colon\Delta\to
    Z$ the closed embedding.  Note that $\widetilde q_{2}\circ i$
    gives an identification $\Delta\simeq M$.  Applying Verdier
    adjunction formula thrice, we get the commutative diagram $$
    \xymatrix{ \reim{(\widetilde q_{2}\circ i)} \epb{(\widetilde
    q_{2}\circ i)} \omega_{M} \ar@{=}[r] \ar@{=}[d] & \omega_M \\
    \reim{\widetilde q_{2}} \reim{i} \epb{i} \epb{\widetilde q_{2}}
    \omega_{M} \ar[r] & \reim{\widetilde q_{2}} \epb{\widetilde q_{2}}
    \omega_{M} \, .  \ar[u] }$$ In other words, the identity of
    $\omega_{M}$ factorizes through $\reim{\widetilde q_{2}}
    \epb{\widetilde q_{2}} \omega_{M} \simeq \reim{q_{2}}\omega_{Z}$.
    One concludes by using~\cite[Exercise~1.4]{Kashiwara-Schapira90}.
\end{proof}

Finally, let us list some functorial properties that the micro-support
enjoys, referring to~\cite{Kashiwara-Schapira90} for proofs.

Consider the correspondence of cotangent bundles associated to $f$:
$$
T^*N \from[\ft] N\times_M T^*M \to[\fp] T^*M.
$$
Let $F\in\BDC(k_{M})$ and assume $f$ is smooth, then $\epb f
F\simeq\omega_{N/M}\tens\opb f F$ and
\begin{equation}
    SS(\opb f F)\subset\ft\fp^\inv \bigl( SS(F) \bigr).
    \label{eq:SSopb}
\end{equation}
Let $G\in\BDC(k_{N})$ and assume $f$ is proper on $\supp(G)$, then
$\reim f G\simeq\roim f G$ and
\begin{equation}
    SS(\roim f G)\subset\fp\ft^\inv \bigl( SS(G) \bigr).
    \label{eq:SSoim}
\end{equation}
Let $F,G\in\BDC(k_{M})$ and assume $SS(F)\cap SS(G)^{a}\subset M$,
then
\begin{equation}
    SS(F\tens G)\subset SS(F) + SS(G).
    \label{eq:SStens}
\end{equation}
Let $F,G\in\BDC(k_{M})$ and assume $SS(F)\cap SS(G)\subset M$, then
\begin{equation}
    SS \bigl( \rhom(G,F) \bigr) \subset SS(F) + SS(G)^a.
    \label{eq:SShom}
\end{equation}

\subsection{Review on kernels}

Consider the natural projections $M\from[q_{1}]M\times M\to[q_{2}]M$.
To $K\in\BDC(k_{M\times M})$ one associates the functors $$
\Phi_{K}(F)=\reim{q_{2}}(K\tens\opb{q_{1}}F), \qquad
\Psi_{K}(F)=\roim{q_{1}}\rhom(K,\epb{q_{2}}F).  $$
These two functors are adjoint to each other, i.e., for
$F,G\in\BDC(k_{M})$
\begin{equation}
    \RHom(\Phi_{K}(F), G) \simeq \RHom \bigl( F, \Psi_{K}(G) \bigr).
    \label{eq:adj}
\end{equation}
Using the estimates we recalled in the previous section, one easily
gets the following result.

\begin{proposition}
    \label{pr:estim}
    Let $F\in\BDC(k_{M})$ and $K\in\BDC(k_{M\times M})$.
    \begin{itemize}
	\item [(i)] Assume that $q_{2}$ is proper on $\supp(K)\cap
	q_{1}^\inv \supp(F)$ and that one has the estimate
	$SS(K)^a\cap (SS(F)\times M)\subset M\times M$.  Then one has
	the estimate $$SS \bigl( \Phi_{K}(F) \bigr) \subset
	\{(y;\eta)\colon (x,y;\xi,\eta)\in SS(K) \text{ for some }
	(x;\xi)\in SS(F)^{a}\}.  $$

	\item [(ii)] Assume that $q_{1}$ is proper on $\supp(K)\cap
	q_{2}^\inv \supp(F)$ and that one has the estimate $SS(K)\cap
	\bigl( M \times SS(F) \bigr) \subset M \times M$.  Then one
	has the estimate $$SS \bigl( \Psi_{K}(F) \bigr) \subset
	\{(x;\xi)\colon (x,y;\xi,\eta)\in SS(K)^{a} \text{ for some }
	(y;\eta)\in SS(F)^{a}\}.  $$
    \end{itemize}
\end{proposition}

\subsection{Proof of Theorem~\ref{th:propag0}}

Let us consider the kernel $K=\reim j\omega_{Z}$, where $j$ denotes 
the embedding $Z\subset M\times M$.  Since $K \simeq 
\rhom(k_{Z},\omega_{M\times M})$, by \eqref{hy:Gconvex} one has $SS(K) 
\subset SS(k_{Z})^a$.  By \eqref{hy:ZtransN} and the fact that $q_{1}$ 
is proper on $Z\cap q_{2}^\inv W$, the hypotheses of 
Proposition~\ref{pr:estim}~(ii) are satisfied.  We find $$SS \bigl( 
\Psi_{K}(F) \bigr) \subset \{(x;\xi)\colon (x,y;\xi,\eta)\in SS(k_{Z}) 
\text{ for some } (y;\eta)\in SS(F)^a\}.$$ Let $(x,y;\xi,\eta)\in 
SS(k_{Z})$ with $(y;\eta)\in SS(F)^a$.  Hypothesis \eqref{hy:ZsubGG} 
together with the fact that $SS(F)\cap\dcone^{a}\subset M$, imply that 
$(y;\eta)\in M$, and then hypothesis \eqref{hy:ZtransM} implies 
$(x;\xi)\in M$.  We thus have $SS \bigl( \Psi_{K}(F) \bigr) \subset 
M$, and hence $\Psi_{K}(F)$ is locally constant on $M$.  On the other 
hand, one has the estimate $\supp \bigl( \Psi_{K}(F) \bigr) \subset 
q_{1}(Z\cap q_{2}^\inv W)=\pas W$, and $\pas W$ does not contain any 
connected component of $M$ by hypothesis.  Hence $\Psi_{K}(F)=0$.

By the same argument we obtain $\Psi_{K}(F\tens\omega_{M})=0$, and
hence $$0 = \RHom(k_{M}, \Psi_{K} \bigl( F\tens\omega_{M}) \bigr)
\simeq \RHom(\Phi_{K}(k_{M}), F\tens\omega_{M}).$$ Since
$\Phi_{K}(k_{M})\simeq\reim{q_{2}}\omega_{Z}$,
Lemma~\ref{le:summand}~(ii) implies $$0 = \RHom(\omega_{M},
F\tens\omega_{M}) \simeq \RHom(k_{M}, F).$$

\begin{remark}
    As it is clear from the above proof, one could generalize the
    notion of propagator by considering pairs $(K,\dcone)$, for
    $K\in\BDC(k_{M\times M})$.  In this case, one should replace
    $Z$-proper by $\supp(K)$-proper, and hypothesis \eqref{hy:Zdelta}
    by the following requirement: there exist $G\in\BDC(k_{M})$ and a
    locally free sheaf of rank one $L$ on $M$, such that $L$ is a
    direct summand of $\Phi_{K}(G)$.
\end{remark}

\subsection{Proof of Corollary~\ref{co:propag}}

Let us prove (i).  Since $SS(F)\cap SS(k_{W})\subset M$, we get by
\eqref{eq:SShom} that $SS(\rsect_{W}F)\subset SS(F)+\dcone$.  Since
$SS(F)\cap\dcone^{a}\subset M$ and $\dcone$ is a proper convex closed
cone, this implies
\begin{equation}
    SS(\rsect_{W}F)\cap\dcone^{a}\subset M.
    \label{eq:SSWF}
\end{equation}
We may then apply Theorem~\ref{th:propag0} with $F$ replaced by
$\rsect_{W}F$.

The proof of (ii) is almost the same, noticing that since $SS(F)\cap
SS(k_{U})^a\subset M$, we get by \eqref{eq:SStens} that
$SS(F_{U})\subset SS(F)+\dcone$.

\subsection{Proof of Corollary~\ref{co:PropagFinite}}

Applying the functor $\rsect(M;\cdot\tens F)$ to the exact sequence $$
0 \to k_{M\setminus N} \to k_{M} \to k_{N} \to 0, $$ we are reduced to
prove $$ \rsect(M; F_{M\setminus N})=0.  $$ By the hypotheses, one has
the isomorphism in $\BDC(k_{M})$
\begin{equation}
    k_{M\setminus N} \simeq \left( 0 \to k_{U_I} \to \cdots \to
    \DSum_{J\subset I,\ |J|=2} k_{U_{J}} \to \DSum_{i\in I} k_{U_{i}}
    \to 0 \right),
    \label{eq:NinUJ}
\end{equation}
where $\DSum_{i\in I} k_{U_{i}}$ is in degree zero.  Hence, it is
enough to prove that $$
\rsect(M;F_{U_{J}})=0 \quad \text{for any } J\subset I. $$ This
follows from Corollary~\ref{co:propag}~(ii).

\section{Applications to hyperbolic systems}

In this section, $M$ is a real analytic manifold, and $k=\C$.

\subsection{Hyperfunction solutions}

We refer to Sato~\cite{Sato59}, Sato-Kawai-Kashiwara~\cite{SKK}, and 
Kashiwara~\cite{Kashiwara95}, for the notions of hyperfunction, 
wave-front set, and $\D$-module, that we shall use.

Let $X$ be a complexification of $M$.
Following~\cite[\S6.2]{Kashiwara-Schapira90}, using the natural
projection $T^*_{M}X\to M$ and the Hamiltonian isomorphism, we will
identify $T^*M$ to a subset of the normal bundle $T_{T^*_{M}X}T^*X$.

Let us denote by $\O_{X}$ and $\D_{X}$ the sheaves of holomorphic 
functions and of linear partial differential operators, respectively.  
If $\shm$ is a coherent $\D_{X}$-module (i.e., a system of PDE), we 
denote by $\chv(\shm)$ its characteristic variety, a closed 
$\C^\times$-conic involutive subvariety of $T^*X$.

\begin{definition}
    (cf~\cite{Kashiwara-Schapira79}) Let $\dcone\subset T^{*}M$ be a
    closed cone, and $\shm$ a coherent $\D_{X}$-module.  One says that
    $\shm$ is {\em $\dcone$-hyperbolic} if $$\dcone \cap C_{T^*_{M}X}
    \bigl( \chv(\shm) \bigr)\subset M.$$
(Note that, since $\chv(\shm)$ is $\C^\times$-conic, $\shm$ is
$\dcone$-hyperbolic if and only if it is $\dcone^a$-hyperbolic.)
\end{definition}

Recall that the sheaf $\shb_{M}$ of Sato's hyperfunctions on $M$ is 
given by $\shb_{M} := \rhom(D'k_{M},\O_{X}) \simeq H^{\dim 
M}_{M}(\O_{X})\tens\ori_{M/X}$.

\begin{theorem}
    \label{th:MhypB}
    Let $(Z,\dcone)$ be a convex propagator on $M$, and $W\subset M$ a
    closed $Z$-proper subset satisfying $SS(k_{W})\subset\dcone^{a}$.
    Let $\shm$ be a coherent $\D_{X}$-module and assume it is
    $\dcone$-hyperbolic.  Then $$\rsect \bigl(
    M;\rhom[\D_{X}](\shm,\sect_{W}\shb_{M}) \bigr) =
    0.$$
\end{theorem}

\begin{proof}
    Setting $F=\rhom[\D_{X}](\shm,\shb_{M})$, it follows
    from~\cite{Kashiwara-Schapira79}
    or~\cite[\S11.5]{Kashiwara-Schapira90} (see
    also~\cite{Bony-Schapira73a} for the case of a single operator)
    that $SS(F)\subset C_{T^*_{M}X} \bigl( \chv(\shm) \bigr)$.  Since
    $\shb_{M}$ is a flabby sheaf, one has $$
    \rsect \bigl( M;\rhom[\D_{X}](\shm,\sect_{W}\shb_{M}) \bigr) =
    \rsect_{W} (M; F).  $$
    The result then follows from Corollary~\ref{co:propag}~(i).
\end{proof}

Let $N\subset M$ be a real analytic submanifold, and denote by
$Y\subset X$ a complexification.  One says that $Y$ is
non-characteristic for a coherent $\D_{X}$-module $\shm$, if
$$
T^*_{Y}X \cap \chv(\shm) \subset X.
$$
In this case, the induced system $\shm_{Y}$ is a coherent
$\D_{Y}$-module.  Note that if $\shm$ is $T^*_{N}M$-hyperbolic, then
$Y$ is non-characteristic for $\shm$.

\begin{theorem}
    \label{th:PropagFiniteB}
    Let $N\subset M$ be a real analytic submanifold, and $I$ a finite
    set.  For $j\in I$, let $U_{j}$ be open subsets of $M$, such that
    $N=M\setminus\Union_{j\in I}U_{j}$.  For any $J \subset I$, $J
    \neq \emptyset$, let $(Z_{J},\dcone_{J})$ be a convex propagator,
    and set $U_{J}=\Inter_{j\in J}U_{j}$.  Assume that $U_{J}$ is
    $Z_{J}$-proper, and $SS(k_{U_{J}}) \subset \dcone_{J}$.  Let
    $\shm$ be a coherent $\D_{X}$-module and assume it is
    $\dcone_{J}$-hyperbolic for any $J \subset I$.  Then, one has the
    isomorphism $$ \rsect \bigl( M;\rhom[\D_{X}](\shm,\shb_{M})
    \bigr)\isoto \rsect \bigl( N;\rhom[\D_{Y}](\shm_{Y},\shb_{N})
    \bigr).  $$
\end{theorem}

Note that the same statement holds when replacing Sato's
hyperfunctions by real analytic functions.

\begin{proof}
    Applying Corollary~\ref{co:PropagFinite} with
    $F=\rhom[\D_{X}](\shm,\shb_{M})$, we get $$
    \rsect \bigl( M;\rhom[\D_{X}](\shm,\shb_{M}) \bigr)\isoto \rsect
    \bigl( N;\rhom[\D_{X}](\shm,\shb_{M})\vert_{N} \bigr).  $$
    It follows by \eqref{eq:NinUJ} that $SS(k_{M\setminus N}) \subset
    \bigcup_{J}\lambda_{J}$.  Since $T^*_{N}M$ coincides with
    $SS(k_{M\setminus N})$ outside of the zero section, the fact that
    $\shm$ is $\dcone_{J}$-hyperbolic for any $J \subset I$ implies
    that $\shm$ is $T^*_{N}M$-hyperbolic.  It then follows
    from~\cite{Kashiwara-Schapira79}
    or~\cite[\S11.5]{Kashiwara-Schapira90} that $$
    \rhom[\D_{X}](\shm,\shb_{M})\vert_{N} \simeq
    \rhom[\D_{Y}](\shm_{Y},\shb_{N}) .  $$
\end{proof}

Let $P$ be a differential operator on $X$, and denote by $\sigma(P)$ 
its principal symbol, a homogeneous function on $T^{*}X$.  One says 
that $P$ is $\dcone$-hyperbolic if so is the associated $\D$-module 
$\shm=\D_{X}/\D_{X}P$.  If $(z)=(x+iy)$ is a local coordinate system 
in $X$, and $(z;\zeta)=(x+iy;\xi+i\eta)$ the associated symplectic 
coordinates in $T^*X$, then $P$ is $\dcone$-hyperbolic if and only if 
$$\sigma(P)(x;i\eta+\theta)\neq 0\qquad \text{for any } (x;i\eta)\in
T^*_{M}X,\ (x;\theta)\in\dcone,\ \theta\neq 0.$$

\begin{corollary}
    \label{co:PhypB}
    Let $(Z,\dcone)$ be a convex propagator on $M$, and $W\subset M$ a
    closed $Z$-proper subset satisfying $SS(k_{W})\subset\dcone^{a}$.
    Let $P$ be a differential operator on $X$ and assume it is
    $\dcone$-hyperbolic.  Then $P$ induces an isomorphism $$P\colon
    \sect_{W}(M;\shb_{M}) \isoto
    \sect_{W}(M;\shb_{M}).$$
\end{corollary}

\begin{proof}
    Apply Theorem~\ref{th:MhypB} with $\shm=\D_{X}/\D_{X}P$, and note
    that the solution complex $\rhom[\D_{X}](\shm,\sect_{W}\shb_{M})$
    is represented by the complex of flabby sheaves
    $$0\to\sect_{W}\shb_{M}\to[P]\sect_{W}\shb_{M}\to 0.$$
\end{proof}

Denote by $r$ the automorphism of $M\times M$ given by $r(x,y)=(y,x)$.

\begin{corollary}
    \label{co:PhypBcauchy}
    Let $N\subset M$ be a real analytic hypersurface dividing $M$ in
    two closed half-spaces $N^\pm$, and let $\theta$ be an analytic
    vector field defined in a neighborhood of $N$ and normal to it.
    Let $(Z,\dcone)$ be a convex propagator on $M$, and assume that
    $N^+$ is $Z$-proper, $N^-$ is $r(Z)$-proper, and
    $SS(k_{N^+})\subset\dcone^{a}$.  Let $P$ be a differential
    operator on $X$, and assume it is $\dcone$-hyperbolic.  Then $P$
    induces a
    surjective morphism $$
    P\colon \sect(M;\shb_{M}) \twoheadrightarrow \sect(M;\shb_{M}), $$
    and moreover the homogeneous Cauchy problem $$
    \begin{cases}
	Pu=0, \\
	\gamma_{\theta}(u)=(w_{1},\dots,w_{m}),
    \end{cases}
    $$
    is globally well posed in the framework of hyperfunctions.  (Here,
    $m$ is the order of $P$, and the trace map $\gamma_{\theta}(u) =
    (u\vert_{N},\theta u\vert_{N},\dots,\theta^{m-1}u\vert_{N})$ is
    well defined since $Pu=0$ implies that the wave-front of $u$ is
    transversal to $N$.)
\end{corollary}

\subsection{Distribution solutions}

As above, let $X$ be a complexification of $M$.  We denote by
$\Db_{M}$ the sheaf of Schwartz distributions on $M$.

\begin{definition}
    Let $\shm$ be a coherent $\D_{X}$-module.
    \begin{itemize}
	\item[(i)] We say that $\shm$ is {\em $\Db$-hyperbolic} at
	$p\in T^*M$ if $$ p \notin SS\bigl(\rhom[\D_{X}](\shm,
	\Db_{M}) \bigr).  $$

	\item[(ii)] Let $\dcone\subset T^{*}M$ be a closed cone.  One
	says that $\shm$ is {\em $\dcone$-$\Db$-hyperbolic} if it is
	$\Db$-hyperbolic at any $p\in\dcone\setminus M$, i.e.\ if
	$$\dcone \cap
	SS\bigl(\rhom[\D_{X}](\shm, \Db_{M}) \bigr) \subset M.$$
    \end{itemize}
\end{definition}

With this definition, it is clear that Corollary~\ref{co:propag}~(i)
implies

\begin{theorem}
    \label{th:MdistrB}
    Let $(Z,\dcone)$ be a convex propagator on $M$, and $W\subset M$ a
    closed $Z$-proper subset satisfying $SS(k_{W})\subset\dcone^{a}$.
    Let $\shm$ be a coherent $\D_{X}$-module, and assume it is
    $\dcone$-$\Db$-hyperbolic.  Then $$\rsect \bigl(
    M;\rhom[\D_{X}](\shm,\rsect_{W}\Db_{M}) \bigr) =
    0.$$
\end{theorem}

\begin{remark}
    The problem, of course, is to give conditions for a system $\shm$
    to be $\Db$-hyperbolic.  If $P$ is a differential operator on $X$,
    and $\shm=\D_{X}/\D_{X}P$, then it is well-known that $\shm$ is
    $\Db$-hyperbolic if: it is hyperbolic, has characteristics with
    real constant multiplicities, and it satisfies the Levi
    conditions.  An analog statement holds for systems (not
    necessarily determined) by~\cite{D'Agnolo-Tonin98}.  Little is
    known beside the case of real constant multiplicities, or of
    constant coefficients in $\R^n$.
\end{remark}

Let us now consider the case of a single differential operator $P$.
One says that $P$ is $\Db$-hyperbolic at $p$
(resp.~$\dcone$-$\Db$-hyperbolic) if so is the system
$\shm=\D_{X}/\D_{X}P$.

\begin{corollary}
    \label{co:PdistrB}
    Let $(Z,\dcone)$ be a convex propagator on $M$, and $W\subset M$ a
    closed $Z$-proper subset satisfying $SS(k_{W})\subset\dcone^{a}$.
    Let $P$ be a $\dcone$-$\Db$-hyperbolic differential operator on
    $X$, and assume that $P\colon \Db_{M}\to \Db_{M}$ is stalk-wise
    surjective.  Then $P$ induces isomorphisms
    \begin{align*}
	P &\colon \sect_{W}(M;\Db_{M}) \isoto \sect_{W}(M;\Db_{M}), \\
	P &\colon H^1_{W}(M;\Db_{M}) \isoto H^1_{W}(M;\Db_{M}).
    \end{align*}
\end{corollary}

\begin{proof}
    Set $\Db_M^P = \ker (P\colon \Db_{M} \to \Db_{M})$.  Since
    $P\colon \Db_{M}\to \Db_{M}$ is an epimorphism, we have an
    isomorphism $\rhom[\D_{X}] (\D_{X}/\D_{X} P ,\Db_{M}) \simeq
    \Db_M^P$, and a short exact sequence $$
    0 \to \Db_M^P \to \Db_{M}\to[P] \Db_{M} \to 0.  $$
    Applying the functor $\rsect_{W}(M,\cdot)$, we get the long exact
    cohomology sequence
    \begin{equation}
	\begin{split}
	0 &\to \sect_W(M;\Db_M^P) \to \sect_W(M;\Db_{M}) \to[P]
	\sect_W(M;\Db_{M}) \\ &\to H^1_W(M;\Db_M^P) \to
	H^1_W(M;\Db_{M}) \to[P] H^1_W(M;\Db_{M}) \\ &\to
	H^2_W(M;\Db_M^P) \to 0.
	\end{split}
	\label{eq:W012}
    \end{equation}
    Theorem~\ref{th:MdistrB} implies $H^j_W(M;\Db_M^P) = 0$ for any
    $j$, and the proof is complete.
\end{proof}

Let us discuss a sufficient condition for $P$ to be $\Db$-hyperbolic.

\begin{proposition}
    \label{pr:Pdistr}
    Let $P$ be a differential operator on $X$, and let $p\in\dot T^*
    M$.  Assume
    \begin{itemize}
	\item[(i)] $\sigma(P)(p) \neq 0$,

	\item[(ii)] $P\colon (\Db_{M})_{x} \to (\Db_{M})_{x}$ is
	surjective for any $x$ in a neighborhood of $\pi(p)$,

	\item[(iii)] there exists an open neighborhood $\Omega\subset
	T^*M$ of $p$ such that for any $x\in\pi(\Omega)$, and any
	$C^\infty$-function $\varphi$ on $M$ with $\varphi(x)=0$,
	$d\varphi(x)\in\Omega$, one has

	(iii)$_{1}$ given $u\in(\sect_{\{\varphi < 0\}}\Db_{M})_{x}$
	satisfying $Pu=0$, there exists $\tilde u\in(\Db_{M})_{x}$
	such that $\tilde u\vert_{\{\varphi < 0\}}=u$ and $P\tilde
	u=0$,

	(iii)$_{2}$ given $v\in(\sect_{\{\varphi < 0\}}\Db_{M})_{x}$
	there exists $u\in(\sect_{\{\varphi < 0\}}\Db_{M})_{x}$ such
	that $Pu=v$.
   \end{itemize}
   Then, $P$ is $\Db$-hyperbolic at $p$.
\end{proposition}

Note that in (i) we used the embedding $T^*M \hookrightarrow M
\times_{M} T^*X$, which exists since $X$ is a complexification of $M$.

\begin{proof}
    Since conditions (i)--(iii) are open in $p\in T^* M$, we may find
    an open neighborhood $\Omega$ of $p$ in $T^*M$ such that (iii)
    holds, (ii) holds in $\pi(\Omega)$, and moreover $\sigma(P)(q)\neq
    0$ for any $q\in\Omega$.  Let $x\in\pi(\Omega)$, and $\varphi$ be
    a $C^\infty$-function on $M$ as in (iii).

    Consider the morphism of exact sequences, where the vertical
    arrows are induced by $P$
    \begin{equation}
	\xymatrix{ (\sect_{\{\varphi \geq 0\}}\Db_{M})_{x}
	\ar@{^{(}->}[r]^-{f} \ar[d]^{\alpha} & (\Db_{M})_{x}
	\ar[r]^-{g} \ar[d]^{\beta} & (\sect_{\{\varphi <
	0\}}\Db_{M})_{x} \ar@{->>}[r]^-{h} \ar[d]^{\gamma} &
	(H^{1}_{\{\varphi \geq 0\}}\Db_{M})_{x} \ar[d]^{\delta} \\
	(\sect_{\{\varphi \geq 0\}}\Db_{M})_{x} \ar@{^{(}->}[r]^-{f} &
	(\Db_{M})_{x} \ar[r]^-{g} & (\sect_{\{\varphi <
	0\}}\Db_{M})_{x} \ar@{->>}[r]^-{h} & (H^{1}_{\{\varphi \geq
	0\}}\Db_{M})_{x} \, .  }
	\label{eq:ABCD}
    \end{equation}
    Consider the stalk-wise analog of \eqref{eq:W012} for $W=\{\varphi
    \geq 0\}$.  By definition of the micro-support, we are left to
    prove that $\alpha$ and $\delta$ are isomorphisms.  This follows
    from the following considerations.  Hypothesis (i) states that
    $\{\varphi = 0\} \subset M$ is non-characteristic for $P$, and by
    Holmgren's theorem this implies that $\alpha$ is injective.  By
    (ii), $\beta$ is surjective.  Moreover, hypothesis (iii)$_{2}$
    says that $\gamma$ is surjective, while hypothesis (iii)$_{1}$
    reads $g\colon \ker\beta \twoheadrightarrow \ker\gamma$.
\end{proof}

\begin{remark}
    In his beautiful paper~\cite{Leray55}, Jean Leray discusses, among 
    other topics, the problem of global extension for solutions to 
    hyperbolic operators with simple characteristics.  In particular, 
    in loc.~cit.\ it is shown that such operators satisfy the 
    hypotheses of Proposition~\ref{pr:Pdistr}.
\end{remark}

\section{Causal manifolds}

\subsection{Conal manifolds}
\label{se:conal}

In this section, we shall construct convex propagators.

\begin{definition}
    We say that a cone $\cone\subset TM$ is {\em admissible} if it is
    closed proper convex and $\Int(\cone_{x})\neq\emptyset$ for any
    $x\in M$. (Here, $\Int(\cone_{x})$ denotes the interior of
    $\cone_{x}$.)  If $\cone \subset TM$ is an admissible cone, we say
    that a closed subset $Z\subset M\times M$ is a {\em
    $\cone$-propagator} if
    \begin{itemize}
	\refstepcounter{equation}\label{hy:Zdelta2} \item
	[\eqref{hy:Zdelta2}] $\Delta\subset Z$,

	\refstepcounter{equation}\label{hy:NZsupGG} \item
	[\eqref{hy:NZsupGG}] $N(Z)\supset \bigl( M \times \Int(\cone)
	\bigr) \cup \bigl( \Int(\cone)^a \times M \bigr) $.
    \end{itemize}
    (As for \eqref{hy:NZsupGG}, recall that we identify the
    zero-section of $TM$ to $M$.)
\end{definition}

\begin{proposition}
    If $\cone \subset TM$ is an admissible cone and $Z \subset M
    \times M$ is a $\cone$-propagator, then $(Z,\cone^{\circ})$ is a
    convex propagator.
\end{proposition}

\begin{proof}
    If $\cone\subset TM$ is admissible, then $\dcone=\cone^{\circ}$
    satisfies \eqref{hy:Gconvex}.

    If $V_{1}$ and $V_{2}$ are two real finite dimensional vector
    spaces, we identify $(V_{1}\times V_{2})^*$ to $V_{1}^*\times
    V_{2}^*$ by $\langle (v_{1},v_{2}), (\nu_{1},\nu_{2}) \rangle =
    \langle v_{1}, \nu_{1} \rangle + \langle v_{2},\nu_{2} \rangle$.
    Then, if $\Vcone_{1}$ and $\Vcone_{2}$ are two cones with
    $\Vcone_{1} \neq \emptyset$, $\Vcone_{2} \neq \emptyset$, one has
    $(\Vcone_{1}\times \Vcone_{2})^\circ = \Vcone_{1}^\circ \times
    \Vcone_{2}^\circ$.  In particular, since
    $\Int(\cone_{x})\neq\emptyset$ for any $x\in M$, one has
    $(M\times\Int(\cone))^\circ = T^*M \times \cone^\circ$.  This last
    set contains $N(Z)^\circ$ by hypothesis \eqref{hy:NZsupGG}.  Using
    the estimate \eqref{eq:SSinN}, \eqref{hy:ZsubGG} follows.

    Remark that if $\Vcone$ is an open convex cone in $V_{1}\times
    V_{2}$ and $(0,v_{2})\in\Vcone$ for $v_{2}\neq 0$, then
    $\Vcone^\circ\cap(V_{1}^*\times\{0\})=\{0\}$.  By hypothesis
    \eqref{hy:NZsupGG}, for each $x_{\circ}, y_{\circ}\in M$ there
    exists $0 \neq w_{\circ}\in T_{y_{\circ}}M$ with $(0,w_{\circ})\in
    N(Z)_{(x_{\circ}, y_{\circ})}$.  Then $N(Z)^\circ \cap
    (T^*_{x_{\circ}}M \times \{y_{\circ}\}) \subset \{0\}$, and
    \eqref{hy:ZtransM} follows.

    The proof of \eqref{hy:ZtransN} is similar.
\end{proof}

\begin{definition}
    A conal manifold is a $C^\infty$-manifold $M$ endowed with an
    admissible cone $\cone \subset TM$.  On a conal manifold $M$, a
    continuous piecewise smooth curve $\alpha\colon [0,1]\to M$ is
    called a {\em $\cone$-path} if the derivative from the right
    $\alpha'_{r}(t)$ exists for any $t\in [0,1[$, and moreover
    $\alpha'_{r}(t) \in \cone_{\alpha(t)}$.  For $x,y\in M$, we write
    $x\preccurlyeq y$ if there exists a $\cone$-path $\alpha\colon
    [0,1]\to M$ with $\alpha(0)=x$, $\alpha(1)=y$.
\end{definition}

Clearly, $\preccurlyeq$ is a preorder relation.  In general, the
graph of $\preccurlyeq$ in $M\times M$ is not closed, and we consider
its closure
\begin{equation}
    Z_{\cone} = \overline{ \{(x,y)\colon x\preccurlyeq y\} }.
    \label{eq:Zcone}
\end{equation}
Note that $Z_{\cone}$ may fail to be the graph of a preordering,
since transitivity may not hold.

\begin{proposition}
    \label{pr:ZconeTrans}
    Let $\cone \subset TM$ be an admissible cone. Then
    \begin{itemize}
	\item [(i)] if $W\subset M$ is a closed subset with $\pas W =
	W$, then $SS(k_{W}) \subset \cone^{\circ a}$.

	\item [(ii)] $Z_{\cone}$ is a $\gamma$-propagator,
    \end{itemize}

\end{proposition}

\begin{proof}
    (i) By \eqref{eq:SSinN}, it is enough to show that $N(W) \supset 
    \Int(\cone^a)$.  Let $x_{\circ}\in W$, and 
    $-v_{\circ}\in\Int(\cone_{x_{\circ}})$.  There exist a local chart 
    $U$ at $x_{\circ}$, and an open conic neighborhood $\Vcone$ of 
    $v_{\circ}$ in $T_{x_{\circ}}M$, such that $U\times \Vcone \subset 
    \cone$.  In view of \eqref{eq:strict}, we shall prove that $$U 
    \cap \bigl( (W\cap U)-\Vcone \bigr) \subset W.$$ Since $W = \pas 
    W$, if $\alpha$ is a $\cone$-path and $\alpha(1)\in W$, then 
    $\alpha(0)\in W$.  Let $x\in W\cap U$ and $v\in \Vcone$ with 
    $x-v\in U$.  Since the segment of straight line from $x-v$ to $x$ 
    is a $\cone$-path, $x-v \in W$.

    (ii) Let us prove that $N(Z_{\cone})\supset M \times \Int(\cone)$.  
    Let $(x_{\circ},y_{\circ}) \in Z_{\cone}$, and $w_{\circ} \in 
    \Int(\cone_{y_{\circ}})$.  Take a local chart $V$ at $y_{\circ}$ 
    and an open conic neighborhood $\Vcone$ of $w_{\circ}$ in 
    $T_{y_{\circ}}M$, such that $V\times \Vcone \subset \cone$.  Let 
    $U\subset M$ be an open neighborhood of $x_{\circ}$.  By 
    \eqref{eq:strict}, for any $x\in U$, $y\in V$, and $w\in \Vcone$, 
    with $x\in \pas y$, and $y+w\in V$, we have to show that $x \in 
    \pas{(y+w)}$.  By definition, $x\in \pas y$ if and only if there 
    exist sequences $x_{n} \to x$, $y_{n} \to y$, with $x_{n} 
    \preccurlyeq y_{n}$ (i.e., there is a $\cone$-path from $x_{n}$ to 
    $y_{n}$).  We may assume $x_{n} \in U$, $y_{n}, y_{n}+w \in V$.  
    Since $w\in \Vcone$, the segment of straight line from $y_{n}$ to 
    $y_{n}+w$ is a $\cone$-path.  Composing the $\cone$-paths above, 
    we get $x_{n} \preccurlyeq y_{n}+w$, which implies $x \in 
    \pas{(y+w)}$ as requested.

    The proof that $N(Z_{\cone})\supset \Int(\cone^a) \times M$ is
    similar.
\end{proof}

\subsection{Causal manifolds}

Recall that we denote by $\Delta\subset M\times M$ the diagonal, and 
by $q_{1}$ and $q_{2}$ the first and second projection from $M\times 
M$ to $M$.  Moreover, for $i,j\in\{1,2,3\}$ let us denote by $q_{ij}$ 
the projection from $M\times M\times M$ to the corresponding factor 
$M\times M$ (e.g., $q_{13}(x,y,z)=(x,z)$).  Recall that a preordering 
$\leq$ on $M$ is determined by its graph $Z=\{(x,y)\colon x\leq y\}$, 
which is a subset of $M\times M$ satisfying
\begin{itemize}
    \refstepcounter{equation}\label{hy:Zreflex} \item
    [\eqref{hy:Zreflex}] $\Delta\subset Z$ (reflexivity),

    \refstepcounter{equation}\label{hy:Ztrans} \item
    [\eqref{hy:Ztrans}] $q_{12}^\inv Z\cap q_{23}^\inv Z\subset
    q_{13}^\inv Z$ (transitivity).
\end{itemize}
One says that $Z$ is a {\em proper preordering} if it is a preordering 
satisfying
\begin{itemize}
    \refstepcounter{equation}\label{hy:Zproper} \item
    [\eqref{hy:Zproper}] $Z\subset M\times M$ is closed, and $q_{13}$
    is proper on $q_{12}^{-1}Z\cap q_{23}^{-1}Z$.
\end{itemize}
Note that the last condition in \eqref{hy:Zproper} means that $\fut
D\cap \pas E$ is compact for any compact subsets $D$ and $E$ of $X$.
In particular, this implies that the intervals $\fut x\cap \pas y$ are
compact.

\begin{definition}
    A {\em causal manifold} $M$ is the data of a manifold $M$ and of
    an admissible cone $\cone \subset TM$ such that the set
    $Z_{\cone}$ in \eqref{eq:Zcone} is a preordering.  If
    moreover \eqref{eq:Zcone} is a proper preordering, $M$ is
    called  a {\em properly causal manifold}.
\end{definition}

\begin{corollary}
    Let $M$ be a properly causal manifold, and $W$ a compact subset of
    $M$.  If $\pas W$ does not contain any connected component of $M$,
    then $\pas W$ is $Z_{\cone}$-proper and satisfies $SS(k_{\pas
    W})\subset\cone^{\circ a}$.
\end{corollary}

In other words, we are in a position to apply
Corollary~\ref{co:propag}.

\begin{proof}
    Hypothesis \eqref{hy:Ztrans} implies that $\pas{(\pas W)} = \pas
    W$.  Hypothesis \eqref{hy:Zproper} implies that $\fut D \cap \pas
    W$ is compact for any compact subset $D$ of $M$.  Finally,
    $SS(k_{\pas W})\subset\cone^{\circ a}$ by
    Proposition~\ref{pr:ZconeTrans}~(ii).
\end{proof}

\subsection{Causal homogeneous spaces}

The toy model for admissible cones is the one considered
in~\cite{Kashiwara-Schapira90}, where $M$ is an open subset of a
vector space $V$, and $\cone = M \times \Vcone \subset TM \simeq
M\times V$ for a closed proper convex cone $\Vcone \subset V$.  In
other words, $\cone$ is a constant cone field.  In this case, using
the notations of section~\ref{se:conal}, $x \preccurlyeq y$ reads $x-y
\in \Vcone$, and $Z_{\cone} = \{(x,y) \colon x-y \in \Vcone\}$.  This
picture is invariant under the group of translations in $V$.

Less trivial examples are obtained by considering other Lie groups.
Let $M=G/H$ be a homogeneous manifold, where $G$ is a real Lie group,
and $H\subset G$ a closed subgroup.  An admissible cone $\cone \subset
TM$ is called {\em invariant} if $\tau_{g}'(x)(\cone_{x}) = \cone_{y}$
for $y = \tau_{g}(x)$, where $\tau_{g}$ denotes the $G$-action on $M$,
$\tau_{g}(\tilde g H)=g\tilde g H$.  One easily proves
(see~e.g.~\cite[\S2.2]{Hilgert-Olafsson97})

\begin{proposition}
    If $\cone \subset TM$ is an invariant admissible cone, then
    $Z_{\cone}$ is the graph of a preordering.
\end{proposition}

Let us denote by $\leq$ the preordering defined by $Z_{\cone}$.
Clearly, this preordering is an {\em invariant ordering}, in the
sense that for any $g\in G$, one has $\tau_{g}(x) \leq \tau_{g}(y)$
whenever $x \leq y$.

Denote by $\mathfrak{g}$ and $\mathfrak{h}$ the Lie algebras of $G$
and $H$, respectively.  Denote by $e\in M$ the equivalence class of
the identity element of $G$.  Noticing that $T_{e}M =
\mathfrak{g}/\mathfrak{h}$, it is clear that the data of an invariant
admissible cone $\cone \subset TM$ is equivalent to the data of a
closed convex cone $\Vcone \subset \mathfrak{g}$ which is invariant by
the adjoint action of $H$, and satisfies $\Vcone \cap \Vcone^a =
\mathfrak{h}$, $\Vcone + \Vcone^a = \mathfrak{g}$.

\begin{definition}
   A {\em causal homogeneous space} $M=G/H$ is the data of a real Lie
   group $G$, a closed subgroup $H\subset G$, and a cone $\Vcone
   \subset \mathfrak{g}$ satisfying the above properties.  $M$ is
   called a {\em properly causal homogeneous space} if the associated
   causal manifold is properly causal.
\end{definition}

If $(G,H)$ is a symmetric pair, refer to~\cite{Hilgert-Olafsson97} for
a wide family of examples of triples $(G,H,\Vcone)$ inducing a
properly causal homogeneous structure on $M=G/H$.

Let us discuss a possible application of our results.

In~\cite{Faraut96}, Faraut constructs global fundamental solutions to
invariant hyperbolic differential operators in the framework of
distributions.  His method relies on the theory of constant
coefficient hyperbolic operators and the technique of spherical
transforms.  Let us show how our results imply the existence of global
fundamental solutions in the framework of hyperfunctions.

Assume that $(G,H,\Vcone)$ induces a properly causal homogeneous 
structure on $M=G/H$.  Let $P$ be an invariant differential operator 
on $M$ such that $$\sigma(P)(e;i\eta+\theta)\neq 0 \qquad\text{for any 
} \eta\in\mathfrak{g}^*,\ \theta\in C^\circ,\ \theta\neq 0.$$ If $\pas 
e$ does not contains the connected component of $e$, we may apply 
Corollary~\ref{co:PhypB} for $W=\pas e$, and get the existence of a 
fundamental solution $$Pu=\delta_{e}, \qquad u\in\sect_{\pas
e}(M,\shb_{M}).$$

{\small

\providecommand{\bysame}{\leavevmode\hbox to3em{\hrulefill}\thinspace}

}


\begin{thebibliography}{1}

\bibitem{Bony-Schapira73a}
Jean-Michel Bony and Pierre Schapira, \emph{Solutions hyperfonctions
du probl{\`e}me de {Cauchy}}, Lecture Notes in Mathematics, no.  287,
Springer-Verlag, 1973, Proceedings Katata 1971, pp.~82--98.

\bibitem{D'Agnolo-Tonin98}
Andrea D'Agnolo and Francesco Tonin, \emph{Cauchy problem for
hyperbolic ${\mathcal D}$-modules with regular singularities}, Pacific
J. Math.  \textbf{184} (1998), no.~1, 1--22.

\bibitem{Faraut96}
Jacques Faraut, \emph{Op\'erateurs diff\'erentiels invariants
hyperboliques sur un espace sym\'etrique ordonn\'e}, J. Lie Theory
\textbf{6} (1996), no.~2, 271--289.

\bibitem{Hilgert-Olafsson97}
Joachim Hilgert and Gestur {\'O}lafsson, \emph{Causal symmetric
spaces}, Perspectives in Mathematics, vol.~18, Academic Press Inc.,
San Diego, CA, 1997, Geometry and harmonic analysis.

\bibitem{Kashiwara95}
M.~Kashiwara, \emph{Algebraic study of systems of partial differential 
equations} (a translation by A.~D'Agnolo and J.-P.~Schneiders of 
Kashiwara's Master's Thesis, Tokyo University, December 1970), M\'em.  
Soc.  Math.  France (N.S.) (1995), no.~63, xiv+72.

\bibitem{Kashiwara-Schapira79}
Masaki Kashiwara and Pierre Schapira, \emph{Micro-hyperbolic systems},
Acta Math.  \textbf{142} (1979), 1--55.

\bibitem{Kashiwara-Schapira90}
\bysame, \emph{Sheaves on manifolds}, Grundlehren der Mathematischen
Wissenschaften, 292, Springer-Verlag, Berlin, 1990.

\bibitem{Leray55}
Jean Leray, \emph{Hyperbolic differential equations}, The Institute
for Advanced Study, Princeton, N. J., 1953 1955.

\bibitem{Sato59}
Mikio Sato, \emph{Theory of hyperfunctions.  I}, J. Fac.  Sci.  Univ.  
Tokyo.  Sect.  I \textbf{8} (1959), 139--193; \emph{II}, ibid.  
\textbf{8} (1960), 387--437.

\bibitem{SKK}
M.~Sato, T.~Kawai, and M.~Kashiwara, \emph{Microfunctions and 
pseudo-differential equations}, Hyperfunctions and pseudo-differential 
equations (Proc.  Conf., Katata, 1971; dedicated to the memory of 
Andr\'e Martineau), Springer, Berlin, 1973, pp.~265--529.  Lecture 
Notes in Math., Vol.  287.

\end{thebibliography}
\end{document}